\documentclass[12pt, a4paper]{article}
\usepackage[left=2.5cm,top=2.2cm,right=2.2cm,bottom=2.2cm,nohead]{geometry}
\usepackage{url}
\usepackage{enumitem}
\usepackage{amsmath, amsfonts, times, setspace, graphicx, lineno, color}
\usepackage{epstopdf}
\usepackage[T1]{fontenc} 
\usepackage{ulem}
\usepackage{caption}

\usepackage[framemethod=TikZ]{mdframed}
\mdfdefinestyle{MyFrame}{%
    linecolor=black,
    outerlinewidth=0.2pt,
    roundcorner=2pt,
    innerrightmargin=20pt,
    innerleftmargin=20pt,
    backgroundcolor=white}

\definecolor{light-gray}{gray}{0.7}
\definecolor{dark-green}{rgb}{0,0.5,0}

\usepackage{amssymb}

\raggedright

\makeatletter
\def\hlinewd#1{%
\noalign{\ifnum0=`}\fi\hrule \@height #1 %
\futurelet\reserved@a\@xhline} \makeatother
\newcommand{\hboldline}{\hlinewd{1pt}}

\usepackage{lineno}
\usepackage{empheq}
\usepackage{mathtools}

\setlist[enumerate]{topsep=0pt}

\begin{document}

\begin{center}
\Large{\bf{Solving multi-objective optimization problems in conservation with the reference point method}}\\
\end{center}

\vspace{30pt}

\begin{flushleft}
Yann Dujardin$^{1,*}$, Iadine Chad\`{e}s$^{1}$\\

\vspace{50pt}

\noindent $^{1}$ CSIRO, Brisbane QLD, Australia\\
$^{*}$ Corresponding author\\
E-mail: yann.dujardin@csiro.au\\
\end{flushleft}



\pagebreak

\section*{Abstract}


Managing biodiversity extinction crisis requires wise decision-making processes able to account for the limited of resources available. In most decision problems in conservation biology, several conflicting objectives have to be taken into account. Most methods used in conservation either provide suboptimal solutions or use strong assumptions about the decision-maker's preferences.
Our paper reviews some of the existing approaches to solve multi-objective decision problems and presents new multi-objective linear programming formulations of two multi-objective optimization problems in conservation, allowing the use of a reference point approach. Reference point approaches solve multi-objective optimization problems by interactively representing the preferences of the decision-maker with a point in the criteria space, called the reference point. We modelled and solved the following two conservation problems: a dynamic multi-species management problem under uncertainty and a spatial allocation resource management problem. 
Results show that the reference point method outperforms classic methods while illustrating the use of an interactive methodology for solving combinatorial problems with multiple objectives. 
The method is general and can be adapted to a wide range of ecological combinatorial problems.

\textbf{Keywords}: conservation; multi-objective; optimization; reference point; linear programming; combinatorial; uncertainty; multi-species dynamic; spatial resource allocation


\section{Introduction}


In recent years, the benefits of using optimization methods to solve decision problems have been widely acknowledged in conservation biology. For example, optimization methods have been developed to best allocate limited resources to protect threatened species \cite{chades2008stop}, protect interacting species \cite{chades2012setting}, design reserves \cite{wilson2006prioritizing,watts2009marxan}, eradicate invasive species \cite{mehta2007optimal}, restore habitat \cite{possingham2015optimal} or translocate species~\cite{rout2013decide}. 
In behavioral ecology, optimization is used to test evolution by natural selection \cite{houston1988dynamic, venner2006dynamic}. 
The control of disease across meta-populations can also be optimized to ensure fastest recovery~\cite{beyer2012implications}. 
Such optimization methods are needed because decision problems are often combinatorial: the possible decisions we have to choose from are combination of smaller ones, which makes the number of possible decisions too large to attempt an exhaustive approach (one cannot generate every possible decision and compare them).

Additionally, many conservation decision problems involve several conflicting objectives~\cite{williams2011adaptive}. For example, when managing interacting species simultaneously in a complex ecosystem, increasing the abundance of one species can result in the decrease of another~\cite{chades2012setting}. Management cost can also be considered as an additional objective. However, these problems are generally converted to single-objective optimization problems, either considering only one objective or considering an a priori aggregation of the objectives~\cite{nicholson2006objectives, chades2008stop, probert2011managing}, but see~\cite{berbel1995application,soleimani2011some} for some exceptions. In contrast to these single-objective approaches, multi-objective combinatorial optimization aims to solve multi-objective combinatorial decision problems without such reduction.

Here, we show that it is possible to solve classic multi-objective combinatorial optimization conservation problems using a cutting edge approach from multi-objective optimization. The reference point method is an interactive approach that provides optimal solutions while accounting for multiple individual objectives. The preferences of the decision-maker are directly expressed as desired values on each objective. These preferences constitute the components of a reference point. Then, an optimization algorithm calculates the closest possible feasible solution to these preferences. If the computed solution is deemed unsuitable, the decision-maker can update his/her preferences and a new solution is calculated. This process can be repeated iteratively until satisfaction of the decision-maker is reached. This type of method is attractive because it does not need any assumptions about the structure of preferences of the decision-maker, i.e.  preferences can be handled even if they are complicated and do not follow a fixed scheme such as a linear trade-offs.
Additionally, associating the reference point method with linear programming as exact underlying optimization method allows us to provide optimal guarantees on solutions computed.

The reference point method has yet to be used in conservation. In this paper, we present the reference point method after introducing some concepts of multi-objective combinatorial optimization and providing a brief critical review of classic approaches. We then demonstrate the benefits of applying the reference point method to two classic combinatorial problems encountered in conservation: a dynamic multi-species decision problem under uncertainty~\cite{chades2012setting} and a spatial resource allocation problem involving several objectives including biodiversity~\cite{higgins2008multi}. We show that the linear programming reference point method outperforms the current approaches used in conservation for solving such multi-objective problems, in term of both optimality and guidance for the decision-maker. Finally, we show that the formulation of the multi-species dynamic problem can be easily extended to any problem using a Markov decision process (MDP) formalism.
 

\section{Method}
\label{sec:method}

\subsection{Multi-objective combinatorial optimization concepts}
\label{sub:RMOO}

Like any decision problem, a single-objective  decision problem has the following ingredients: a model, a set of controls (called variables), and an objective function depending on the variables~\cite{walters1978ecological}. Additionally, in conservation, and in ecology in general, decision problems may seek to maximize several objectives simultaneously~\cite{walters1978ecological}. It is then worth considering the formal formulation of multi-objective combinatorial optimization problems~\cite{vanderpooten1990multiobjective}:

\begin{align*}
\label{multiobj_pb}
\tag{P}
\max \quad & f_{1}(x),\dots,f_{p}(x)\\
s.t. & \qquad x\in X
\end{align*}

where $f_{j},\, j=1,\ldots,p$, $p\geq2$, are the objectives (or criteria), $x$ is the vector of decision variables which can only take value in the set $X$ of \textit{feasible} (i.e. \textit{possible}) decisions. Because we are in a combinatorial context, $X$ is assumed to be discrete.

Any decision $x\in X$ matches with a point $z\in Z_X = \{(f_{1}(x),\dots,f_{p}(x))~|~x\in X\}$.
In contrast to single-objective optimization problems, which admit at most one optimal value, multi-objective optimization problems often admit several optimal points, i.e. points of $Z_X$ that cannot be outperformed on every objective simultaneously by another point of $Z_X$. These points are called \textit{non-dominated points}.
Formally, a non-dominated point is a point $z\in Z_{X}$ such that there is no $z'\in Z_{X}$ with the property $z'\geq z$, where inequality $\geq$ between two points of $Z$ is defined in Box~1.

\bigskip
\begin{figure}[h!]
\begin{mdframed}[style=MyFrame]
For every $z,\, z'\in Z$, $z'\geq z$ if and only if for all $j\in\left\{ 1,\ldots,p\right\} ,\, z_{j}'\geq z_{j}$ and there is $k\in\left\{ 1,\ldots,p\right\}$ such that $z_{k}'>z_{k}$.
\end{mdframed}
{\caption*{Box 1}}
\end{figure}




Non-dominated points are essential in multi-objective combinatorial optimization, since they represent the best possible compromises between the different criteria. The set of non-dominated points and is also called the \textit{Pareto frontier}~\cite{pareto1896course}. Consequently, a \textit{Pareto frontier} represents the set of the best compromise solutions between the different objectives. Most of the multi-objective combinatorial optimization approaches aim to discover non-dominated points and their corresponding decisions, called \textit{efficient} decisions. Indeed, multi-objective combinatorial optimization is often related to one of the following well-known underlying challenges:
\begin{enumerate}
\item Find a particular non-dominated point of the Pareto frontier, according to the preferences of a decision-maker (called \textit{local} approach in this paper);
\item Discover the entire Pareto frontier, or an approximation of this set (called \textit{global} approach in this paper).
\end{enumerate}

 
When the number of criteria is large (>3), it becomes difficult to calculate, represent and analyze the Pareto frontier. Consequently, the local approach should be preferred for problems where the number of criteria may be more than 3.

Finding a non-dominated point according to preferences of decision-makers can be tackled using an \textit{aggregation function}~\cite{ehrgott2006multicriteria}, sometimes also called a \textit{scalarizing function}. The role of aggregation functions is to discriminate non-dominated points according to some preferences. More precisely, an aggregation function is a function $s$ from $Z$ to $\mathbb{R}$, which associates a unique real value to every point of the criteria space. In multi-objective combinatorial optimization, $s$ also depends on parameters called  \textit{preferential parameters}, representing the preferences of a \textit{decision-maker}~\cite{VDP-Interactive-1989}. 
The decision-maker can be a person, a group of persons or any entity able to provide preferences.

\subsection{Classic multi-objective optimization approaches in conservation}\label{sub:classic_eco}

Multi-objective optimization has been used for a long time in fisheries~\cite{mardle1999review}, forestry~\cite{diaz2008making} and natural resources management~\cite{higgins2008multi}. In these fields multi-objective optimization is referred to \textit{multi-objective programming} if not interactive and \textit{interactive processes} otherwise. \textit{Goal programming},  and \textit{compromise programming}, which aim to minimize the deviation between the achievement of goals and their aspiration levels (fixed by the decision-maker in the goal programming case and computed in the compromise programming case), are also popular in these fields~\cite{mardle1999review,diaz2008making}. 
This section will focus on multi-objective optimization in \textit{conservation}. Multi-objective optimization is less developed in conservation than in forestry or fisheries. For example, the well-known approximate conservation solver Marxan \cite{watts2009marxan} is not a multi-objective solver, because the multiple objectives called "targets" are considered as constraints and not as objectives, and no multi-objective optimization framework is yet considered.

\subsubsection{Explicit approaches}\label{sub:expl_app}

Finding optimal solutions when explicitly accounting for multiple objectives in combinatorial problems is a mathematically challenging endeavor. A way to avoid this mathematical challenge is to use \textit{explicit} approaches, i.e. generate a few feasible solutions and compare their performance either by sampling using a model~\cite{geneletti2008protected,schwenk2012carbon} or empirically by asking experts~\cite{carwardine2012prioritizing,joseph2009optimal}.
Although this approach is not, strictly speaking, multi-objective optimization, it is very common in conservation. 

\paragraph*{MCDA}
~\\~\\

The explicit approach allows us to perform multi-criteria decision analysis (MCDA), which is very powerful where the number of possible decisions is small~\cite{lahdelma2000using, kiker2005application, mendoza2006multi}.
The goal of MCDA methods is to determine a best decision or strategy among a reasonable number of possible ones, given that these decisions/strategies are evaluated on several criteria

\paragraph*{Trade-off analysis}
\label{par:trade-off}
~\\~\\
Another usual approach in conservation is to try to establish correlations between criteria (trade-off analysis), via exhaustive approaches~\cite{pichancourt2014growing}, or heuristic approaches~\cite{venter2009harnessing}. 
Unfortunately, there is no reason that criteria of combinatorial problems have the same correlation from one instance to another (changing the data could result in a complete different correlation).
Additionally, the lack of scalability of exhaustive approaches and the lack of optimality of heuristic approaches make them very limited  approaches to solve combinatorial problems. 

%
%
%
%

\subsubsection{Implicit approaches}
\label{sub:aggreg}

When the multi-objective problem can only be implicitly defined (see Section~\ref{sub:RMOO} for a formal definition), we are then confronted to a multi-objective optimization problem.

Local approaches can be used to perform two types of methods: \textit{a priori} methods and interactive methods. \textit{A priori} methods use a unique aggregation function, fixed and defined once by the decision-maker, while interactive methods allow the decision-maker to iteratively change his/her preferences. Global approaches, which generate the entire Pareto frontier, are often called \textit{a posteriori} methods.




\paragraph*{\textit{A priori} methods}
~\\

Several approaches in conservation aim to find a unique objective summarizing the individual objectives, and then treat the problem as a single objective. Reducing several objectives in one is usually done using an \textit{a priori} aggregation function, i.e. an aggregation function with fixed preference parameters. 
The cost-benefit approach is probably the most used approach applying this principle. The cost-benefit approach is an economic approach where every criteria is considered as having an economic counter-part \cite{beria2012multicriteria}. 
Such functions are often used to perform a "cost-benefit" analysis~\cite{chades2008stop,hauser2009streamlining}. Other aggregation functions of the objectives have been studied in conservation~\cite{nicholson2006objectives, giljohann2015choice}. 
Several major well-known drawbacks occur in these approaches. Using economic values of species is ethically controversial because it requires associating an economic value to species~\cite{spash2015bulldozing}. Additionally, in practice, depending on the economical evaluation methods, the value of a species can vary significantly, sometimes from one to tenfold~\cite{loomis1996economic}.
The second drawback is related to the subjectivity and the complexity of the fixed aggregation function. Choosing among a set of potential aggregation functions can be difficult to justify~\cite{nicholson2006objectives, giljohann2015choice}. Finally, reasoning with one objective (aggregated function) instead of several, reduces considerably the role of the decision-maker in the optimization process. Indeed, his/her role is then limited to define the problem. A prescribed solution is then provided by the scientists, missing an opportunity to involve the decision-maker in the decision-making process itself.


%
%

\paragraph*{\textit{A posteriori} methods}
~\\
\textit{A posteriori} methods aim to generate the Pareto frontier or an approximation of it. Some methods in conservation can be classified as \textit{a posteriori} methods \cite{cheung2008trade,beyer2016solving}. \textit{A posteriori} methods corresponds to trade-off analysis methods for implicit approaches. Generating the Pareto frontier is only possible and relevant for problems with a small number of objectives.

\paragraph*{Interactive methods}
\label{sub:RP_method}
~\\

Generally based on the use of parametric aggregation functions, \textit{interactive methods} aim to find a "best compromise solution" interactively with decision-makers \cite{ehrgott2006multicriteria,VDP-Interactive-1989}. 
In these methods, the decision-maker preferences can evolve according to the following iterative procedure:
\begin{itemize}
\item{} Optimization results are obtained using current preferences;
\item{} New preferences are obtained by eliciting feedback from the decision-maker on current results.
\end{itemize}


Fig 1 provides a graphical illustration of the interactive procedure.

\begin{figure}[h!]
\begin{center}
\includegraphics[width=0.7\textwidth]{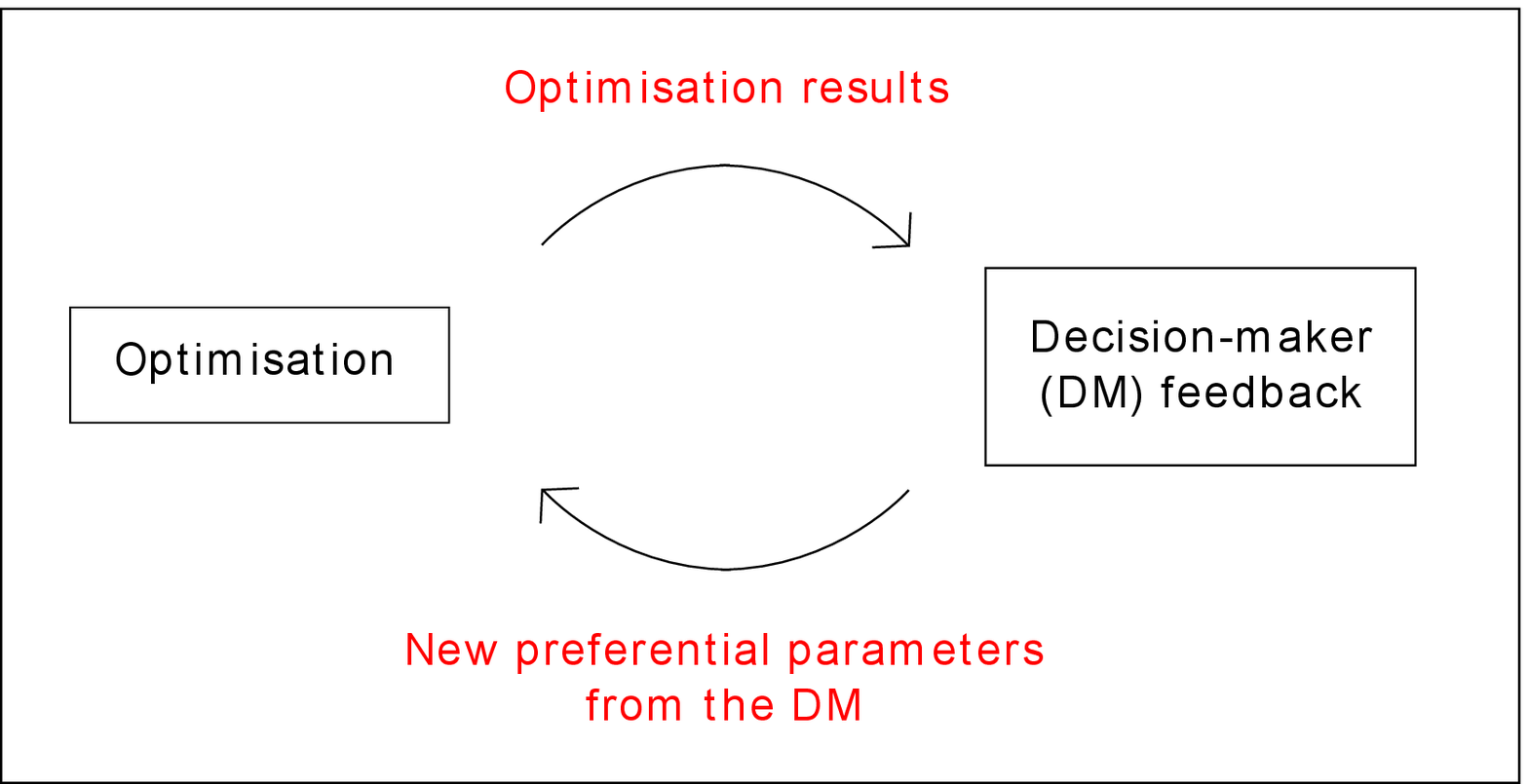}
{\caption*{Fig 1. Multi-objective combinatorial optimization interactive procedure}}
\end{center}
\end{figure}

Multi-objective optimization interactive methods are not very common in conservation but see~\cite{soleimani2011some} for an exception. 

In theory, every \textit{a priori} method can be used in an interactive way. Unfortunately, none of the previously mentioned \textit{a priori} method can satisfy the following requirements: \begin{itemize}
\item Every solution provided by the algorithm corresponds to a non-dominated point of the multi-objective problem.
\item Every non-dominated point of the multi-objective problem can be generated by the algorithm.
\end{itemize}
These requirements are very important, because we would like to avoid generating points which are dominated by others, while being able to generate every possible non-dominated point for the seek of a good compromise solution.

Maybe the most famous interactive method is the weighted sum method. The (interactive) aggregation function is the weighted sum of the normalized objectives where weights represent only the "importance" of every criterion. This method unfortunately fails to satisfy the requirements. In~\cite{probert2011managing}, the authors use a weighted sum method to aggregate two objectives: learning about the model and managing in an optimal adaptive management framework (not used in an interactive way however).
If the criteria are normalized, the weights can be seen as preference parameters defined in Section~\ref{sub:RMOO}. However, several drawbacks still occur:
The first drawback is that weights have no \textit{significance}, and transforming them into meaningful values~\cite{marler2010weighted} can be obscure for the decision-maker~\cite{roy1991outranking}, because the true preferential parameters are hidden.
More importantly, the weighted sum method is well-known to not provide good \textit{compromise solutions}~\cite{pernyRPMDP2010}, i.e. solutions which are well balanced when considering their criteria values. The weighted sum method favors extreme solutions.
Finally, this method makes the strong assumption that one objective can always linearly compensate linearly another objective. However, as raised in \cite{keeney2002common}, the following question has in general no answer: "how much must be gained in the achievement of one objective to compensate for a lesser achievement on a different objective?". Human preferences are often much more complicated than linear trade-offs and may require more elaborate methods.

Conversely, the reference point method is one of the only multi-objective optimization methods to satisfy the requirements~\cite{vanderpooten1990multiobjective}.
Because in this paper we use only the linear programming formulation of the reference point method, we will only present a linear programming formulation (Section~\ref{sub:RPLP}).

\subsection{The linear programming formulation of the reference point method}
\label{sub:RPLP}

Linear programming is well-known in conservation~\cite{chernomor2015split}, but less in the context of multi-objective optimization, especially for reference point based approaches.

Using the reference point method requires solving an optimization problem at every iteration of the interactive process. Although other optimization methods are possible, linear programming (LP) is particularly well suited to solve this problem. 
The LP formulation with $p$ objectives, $n$ variables and $m$ constraints is:
\\
\begin{equation}\label{LP_Multi}
\begin{aligned}
& \max~~~ z + \rho\sum_{j=1}^p \lambda_j (f_j(x) - \overline{z}_j) \\
& s.t.~~~~~~~  z\leq \lambda_j (f_j(x) - \overline{z}_j), j\in \{1,...,p\} \\
& ~~~~~~~~~~~~~a_i\cdot x \leq b_i, i\in \{1,...,m\}\\
& ~~~~~~~~~~~~~x\in \mathbb{Z}^n\\
& ~~~~~~~~~~~~~z\geq 0 
\end{aligned}
\end{equation}
\\
Variables are components of vector $x$, $z$. $f_j, j\in \{1,...,p\}$ are linear functions. $x$ represents the decision, $f_j(x), j\in \{1,...,p\}$ are the objective values, and $z=\min_{j\in \{1,...,p\}}{f_j(x)}$.
Every inequality $a_i\cdot x \leq b_i$ represents a constraint of the problem, while $z\leq \lambda_j (f_j(x) - \overline{z}_j)$ are constraints implying that $z=\min\{\lambda_j(f_j(x)-z_j)\}$.
$\lambda_j = \frac{1}{z_j^{max}-z_j^{min}}$ are fixed and play the role of a normalizing factor.  For each criterion $j$, $z_j^{max}$ and $z_j^{min}$ are respectively the maximum and the minimum possible values of $f_j(x)$ obtained by performing the corresponding single-objective optimization. $\rho$ is a small strictly positive number required to avoid generating weakly non-dominated points, i.e. points which can be dominated over a subset of objectives. Avoiding weakly dominated points is possible by setting $\rho$ to a value inferior to $\frac{\min_{j\in \{1,\cdots,p\}{\lambda_j}}}{\sum_{j\in \{1,\cdots,p\}}{(z_j^{max}-z_j^{min})}}$~\cite{duj2013thesis}, thanks to the combinatorial context ($x$ takes discrete values). 
\\
The simplicity of the formulation is one reason of the popularity of LP to solve the reference point method, which has been implemented in many fields where combinatorial problems occur, for example in telecommunication \cite{ogryczak2005telecommunications}, finance \cite{mansini2014twenty} or transportation \cite{duj2015}. 

\subsection{Adapting the reference point method to two classes of conservation problems}
\label{sub:adapting}
Unfortunately, the LP reference point method requires finding good LP formulations, which in the case of discrete optimization, is in general a hard task and requires a strong knowledge of the conservation problems and LP techniques. 
In this section we present LP formulations for a multi-species dynamic conservation problem and a multi-objective environmental spatial resource allocation problem, so that the reference point method can be applied in both cases.
\\
\subsubsection{Dynamic problem in conservation \label{sub:dyn_cons}}
\label{sub:methDPC}
In \cite{chades2012setting}, the authors propose a method for solving a sequential decision problem under uncertainty, aiming to conserve simultaneously two interacting endangered species: Northern abalone and sea otters. This bi-objective is solved using (indirectly) an \textit{a priori} weighted sum of the objectives. Different weights are tested to generate and explore alternatives. Weighting the criteria allows the use of classic MDP solution methods such as dynamic programming \cite{chades2014mdptoolbox}.
More specifically, the problem is a predator-prey problem where interactions between sea otters and their preferred prey abalone are described using a MDP formalism. Every year, managers must decide between 4 actions: introduce sea otters, enforce abalone anti-poaching measures, control sea otters, half enforce anti-poaching measures and half control sea otters. The time horizon is 20 years. The original problem aims to maximize the density of abalone and abundance of sea otters. 
\\
Because weighting the objectives of an optimization problem can be controversial (see Section~\ref{sub:classic_eco}), we propose to use the linear programming reference point method. 
Adapting Program~\ref{LP_Multi} directly to a LP formulation of MDPs is challenging because rewards appear only in the constraints and not in the objective.
In \cite{pernyRPMDP2010}, the authors were confronted with the same situation when they tried to apply a similar multi-objective optimization technique (the Chebyshev method) to MDPs in a robotic context. The Chebyshev method aims to minimize the Chebyshev norm between the reference point and the decision space. The reference point method is different for several reasons. First, in the reference point method, preferences of the decision-maker are directly expressed as values on every criterion, while in the Chebyshev method preferences are expressed as weights. Second, the reference point method allows the decision-maker to choose values inside the feasible space, which is not the case in the Chebyshev method. However, the same idea as in \cite{pernyRPMDP2010} can also be used for the linear programming reference point method formulation. We first wrote the single-objective dual LP formulation of MDPs and then adapted Program~\ref{LP_Multi} to it. LP and dual LP formulations of MDPs are available in~\cite[chap. 4]{bertsimas1997introduction}.
\\
Formally, the multi-objective Markov decision process related to our multi-species problem is defined by the t-uplet: $\{S,A,H,R_{A},R_{SO},Tr\}$. $S$ is the state space, $A$ is the action space and $H=\{0,\cdots,T-1\}$ is the time-horizon of size $T$. Taking action $a\in A$ when in state $s\in S$ leads to an immediate reward $R_A(s,a)$  \textit{for abalone} and $R_{SO}(s,a)$ \textit{for sea otters}. $Tr$ is the transition matrix. Further details and values are available in \cite{chades2012setting}. Program \ref{LPDP} is the linear programming formulation reference point method we wrote.
\\
\begin{equation}\label{LPDP}\tag{$LP_{DP}$}
\begin{aligned}
& \max~~~ z + \rho (\lambda_A (C_A - \overline{C}_A) + \lambda_{SO} (C_{SO} - \overline{C}_{SO})) \\
& s.t.~~~~~~~  z\leq \lambda_A (C_A - \overline{C}_A) \\
& ~~~~~~~~~~~~~z\leq \lambda_{SO} (C_{SO} - \overline{C}_{SO}) \\
& ~~~~~~~~~~~~~C_A \leq \sum_{t\in H, s\in S, a\in A}{R_A(s,a)x_{t,s,a}} \\
& ~~~~~~~~~~~~~C_{SO} \leq \sum_{t\in H, s\in S, a\in A}{R_{SO}(s,a)x_{t,s,a}} \\
& ~~~~~~~~~~~~~\sum_{a\in A}{x_{t,s,a}} - \sum_{s'\in S, a\in A}{Tr(s',a,s)x_{t,s,a}}=0,~t\in H,~s\in S \\
& ~~~~~~~~~~~~~x_{t,s,a} \geq 0,~t\in H,~s\in S,~a\in A \\
& ~~~~~~~~~~~~~C_A \geq 0\\
& ~~~~~~~~~~~~~C_{SO} \geq 0\\
\end{aligned}
\end{equation}
\\
The main variables are the dual variables $x_{t,a,s}$ of the initial problem. Variables $C_A$ and $C_{SO}$ represent respectively the normalized density of abalone over 20 years and the normalized number of sea otters over 20 years. $(\overline{C}_A, \overline{C}_{SO})$ is the reference point which corresponds to the current preferences of the decision-maker. Note that this LP formulation can easily be generalized for any multi-objective Markov decision process problem, which makes our approach very general (see Box 2).
\\
\begin{empheq}[box=\fbox]{equation}
\begin{aligned}
\nonumber
& \max~~~ z + \rho \sum_{j\in J}\lambda_j (C_j - \overline{C}_j) \\
& s.t.~~~~~~~  z\leq \lambda_j (C_j - \overline{C}_j), j\in J \\
& ~~~~~~~~~~~~~C_j \leq \sum_{t\in H, s\in S, a\in A}{R_j(s,a)x_{t,s,a}}, j\in J \\
& ~~~~~~~~~~~~~\sum_{a\in A}{x_{t,s,a}} - \sum_{s'\in S, a\in A}{Tr(s',a,s)x_{t,s,a}}=0,~t\in H,~s\in S \\
& ~~~~~~~~~~~~~x_{t,s,a} \geq 0,~t\in H,~s\in S,~a\in A \\
& ~~~~~~~~~~~~~C_j \geq 0\\
& ~~~~~~~~~~~~~\text{$J$ is the set of objective}\\
& ~~~~~~~~~~~~~\text{$R_j$ is the reward function associated with objective $j\in J$}\\
& ~~~~~~~~~~~~~\text{$\overline{C}_j$ is the current preference on objective $j\in J$}\\
& ~~~~~~~~~~~~~\text{$Tr$ is the transition matrix of the MDP}\\
\\
\end{aligned}
\end{empheq}
\begin{figure}[h!]
{\caption*{Box 2}}
\end{figure}
\\

\subsubsection{Spatial allocation of resources}
\label{sub:methSAR}

Spatial allocation of resources is an important challenge in conservation including, but not limited to, reserve design~\cite{watts2009marxan,moilanen2008zonation} or environmental investment decision making problems~\cite{higgins2008multi}. In this section, we provide a linear programming reference point formulation of the problem, and demonstrate the use of the reference point method to tackle a spatial resource allocation problem.
\\
In our model, we consider an environmental investment decision making problem inspired by~\cite{higgins2008multi}. We considered a map of 3600 cells, where a decision consists in selecting a subset of 120 cells for management under a budget constraint. In \cite{higgins2008multi} only three objectives were considered, which allows an \textit{a posteriori} approach. As discussed in Section~\ref{sub:RMOO}, this approach has limitations. In particular, a posteriori approaches are relevant only for few criteria (e.g. 3 criteria) while the reference point approach can deal with a large amount of criteria.
\\
We extended the model proposed in \cite{higgins2008multi} by considering five criteria.
The first criterion is related to the minimisation of the total travel time of water. Selected cells will benefit from  management allowing prevention of fast runoff from the highest cells to the water points. The second criterion is related to the maximisation of carbon sequestration. Every selected cell contributes to an improved carbon sequestration in different ways. These two criteria are explained in details in~\cite{higgins2008multi}. We define three additional criteria related to biodiversity. Each of these criteria represent the contribution of the selected cells to the conservation of a different species. 
\\
We considered a map of $|I|\times |J|$ cells where $I=J=\{1,...,60\}$.
We first generated an elevation map, i.e. for every cell $(i,j)\in I\times J$ we generated an elevation $e_{i,j}$. 
According to this elevation map, water runoffs were computed, such that for every cell, the water comes from the highest neighbor (in case of several highest neighbors, one is picked randomly). Thus, every cell $(i,j)$ has a unique antecedent $A((i,j))$, except the peaks of the map which have no antecedent, where we set $A((i,j))=\emptyset$. 
\\
\paragraph*{Main variables}
~\\
For every cell $(i,j)$, $x_{i,j}$ is a 0-1 variable taking the value 1 if $(i,j)$ is managed, and the value 0 otherwise.
\paragraph*{Water Traveling Time criterion}
~\\
For every cell $(i,j)$, $t_{i,j}$ is the average time the water stays on the cell when not managed. $d_{i,j}$ is the additional time water stays on the cell when managed. In our experiments, $t_{i,j}$ and $d_{i,j}$ are random values. For every cell $(i,j)$, $T_{i,j}$ is the time for water to travel the path from the origin cell to the cell $(i,j)$ . $T_{i,j}$ is then equal to the time needed to reach the antecedent cell $A((i,j))$ plus the time of staying on the cell. The Water Traveling Time criterion $WTT$ is the total time needed for water to reach every cell. 
\paragraph*{Carbon sequestration criterion}
~\\
For every cell $(i,j)$, managing $(i,j)$ increases its carbon sequestration value by $c_{i,j}$. The carbon sequestration criterion $CS$ is equal to the sum of $c_{i,j}$ over the managed cells $(i,j)$. 
\paragraph*{Biodiversity criteria}
~\\
For every cell $(i,j)$ and every species $S\in \{1,2,3\}$, managing $(i,j)$ increases the number of individuals of species $S$ by $n_{i,j}^S$. For every species $S$, the biodiversity criteria $N_S$ is equal to the total number of the saved individuals by management. 
\paragraph*{Budget constraint}
~\\
Finally, the cost of managing any cell $(i,j)$ is denoted by $cost_{i,j}$. The management is constrained to respect a budget $B$.
\paragraph*{Linear program}
~\\
\ref{LPRA} below is the linear program associated to the multi-objective resource allocation problem considering the 5 criteria $WTT$, $CS$ and $N_S,~S\in \{1,2,3\}$, and a budget equal to $B$.~\ref{LPRA} is the application of Program~\ref{LP_Multi} to our resource allocation problem. $WTT$ is represented by variable $z_{WTT}$. $CS$ is represented by variable $z_{CS}$. Each $N_S$ is represented by variable $z_{N_S}$.
\\
\begin{equation}\label{LPRA}\tag{$LP_{RA}$}
\begin{aligned}
& \max~~~ z + \rho [ (\lambda_{WTT} (z_{WTT} - \overline{z}_{WTT}) + \lambda_{CS} (z_{CS} - \overline{z}_{CS})+\sum_{S\in \{1,2,3\}}\lambda_{N_s} (z_{N_S} - \overline{z}_{N_S})) ] \\
& s.t.~~~~~~~  z\leq \lambda_{WTT} (z_{WTT} - \overline{z}_{WTT}) \\
& ~~~~~~~~~~~~~z\leq \lambda_{CS} (z_{CS} - \overline{z}_{CS}) \\
& ~~~~~~~~~~~~~z\leq \lambda_S (z_S - \overline{z}_S), S\in \{1,2,3\} \\
& ~~~~~~~~~~~~~\sum_{(i,j)\in I\times J} cost_{i,j}x_{i,j} \leq B \\
& ~~~~~~~~~~~~~z_{WTT} \leq \sum_{(i,j)\in I\times J} T_{i,j}\\
& ~~~~~~~~~~~~~T_{i,j} \leq T_{A((i,j)} + t_{i,j}+ x_{i,j}d_{i,j}, (i,j)\in I\times J\\
& ~~~~~~~~~~~~~z_{CS} \leq \sum_{(i,j)\in I\times J} c_{i,j}x_{i,j} \\
& ~~~~~~~~~~~~~z_{N_S} \leq \sum_{(i,j)\in I\times J} n_{i,j}^Sx_{i,j},~S\in \{1,2,3\} \\
& ~~~~~~~~~~~~~x_{i,j}\in \{0,1\}, (i,j)\in I\times J\\
& ~~~~~~~~~~~~~z_{WTT} \geq 0\\
& ~~~~~~~~~~~~~z_{CS} \geq 0\\
& ~~~~~~~~~~~~~z_{N_S} \geq 0, S\in \{1,2,3\}\\
\end{aligned}
\end{equation}
\\
Our approach is exact and accounts for more objectives than in~\cite{higgins2008multi}. One can also compare the optimal solutions with the solutions provided by the usual explicit approaches. Given the combinatorial nature of the problem, an exhaustive search is of course not possible. 
We tested a possible explicit approach consisting in randomly generating 10,000 decisions respecting the budget constraint. From these decisions we kept 300 points which are non dominated by other generated points. In doing so, our aim is perform a MCDA approach using the 300 points as possible decisions. 
\\
We applied the reference point method using every generated point of the explicit approach as a reference point. In other words, we projected the points of the explicit approach to the Pareto frontier using our LP formulation.

\section{Results}\label{sec:results}


For both case studies, we used the optimization solver Cplex (version 12) to solve the corresponding linear programs.

\subsection{Dynamic problem in conservation \label{sub:dyn_cons}}
\label{sub:resDPC}

Our experiments consisted in comparing the weighted sum approach to the reference point approach (see Section~\ref{sub:methDPC}). 

Fig 2(a) shows the resulting non-dominated points using the weighted sum method applied to our dynamic conservation problem. Twenty equally distributed pairs of weights from $(0,1)$ to $(1,0)$ were used, generating only 4 distinct non-dominated points (the 20 points match to the 4 distinct points). In the context of an interactive procedure, the guidance provided to the decision-maker is then limited. Additionally, none of the non-dominated points represents a good compromise solution between the two objectives since no point has similar values on x-axis and y-axis. 


Fig 2 (b) shows the resulting non-dominated points using the reference point method applied to the same problem and using 20 equally distributed reference points in the criteria space: we computed the extreme points A and B of the Pareto frontier and  subdivided the segment [AB] into 20 points. This time 19 distinct non-dominated points were obtained and interesting good compromise solutions can be identified (similar values on both criteria). In the context of an interactive procedure, the guidance provided to the decision-maker offers a higher chance of reaching satisfaction thanks to the number of distinct points.

\begin{figure}[h!]
    \centering
    \begin{tabular}{cc}
	\includegraphics[width=0.5\textwidth]{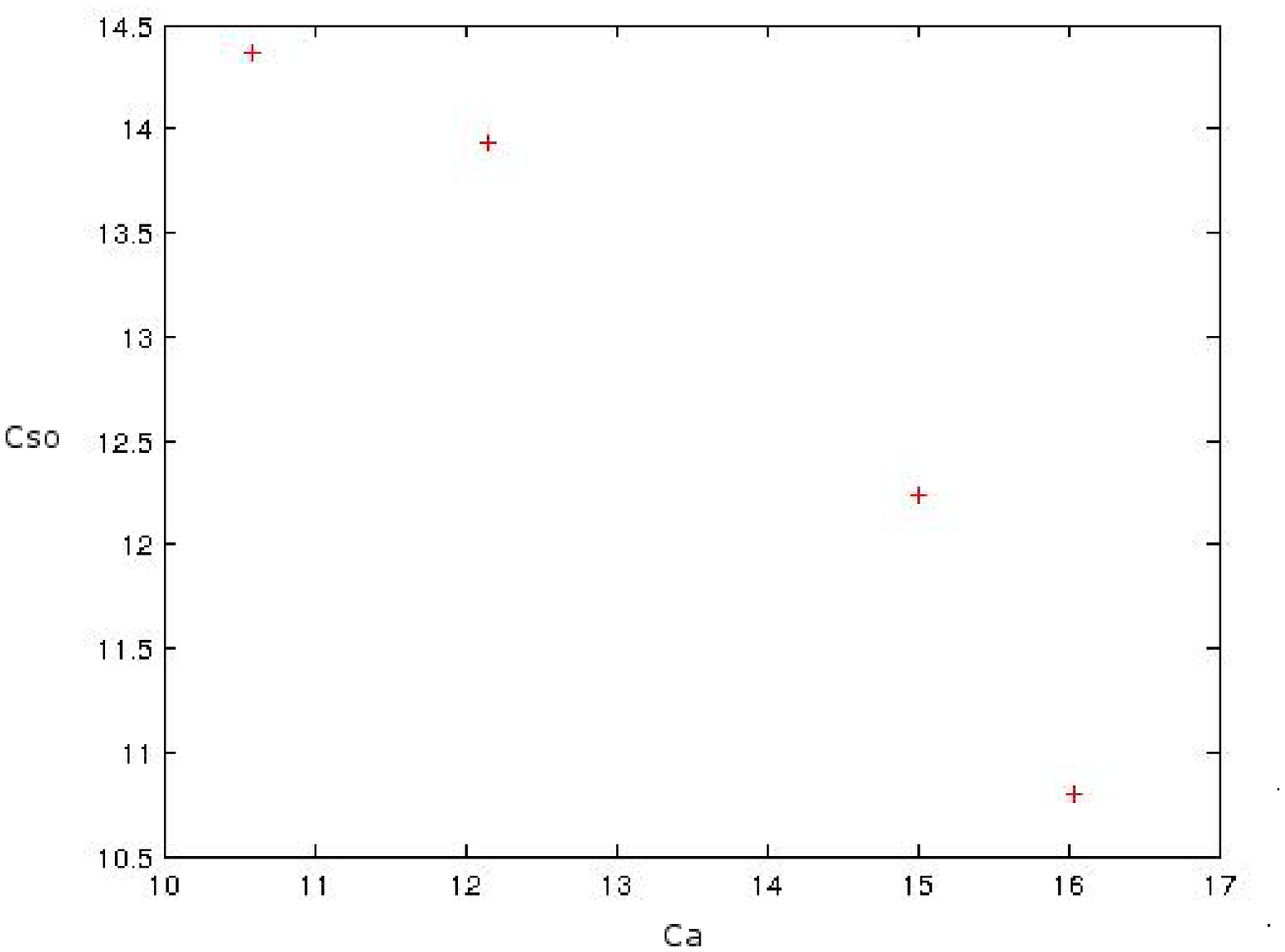} & 
	  \includegraphics[width=0.5\textwidth]{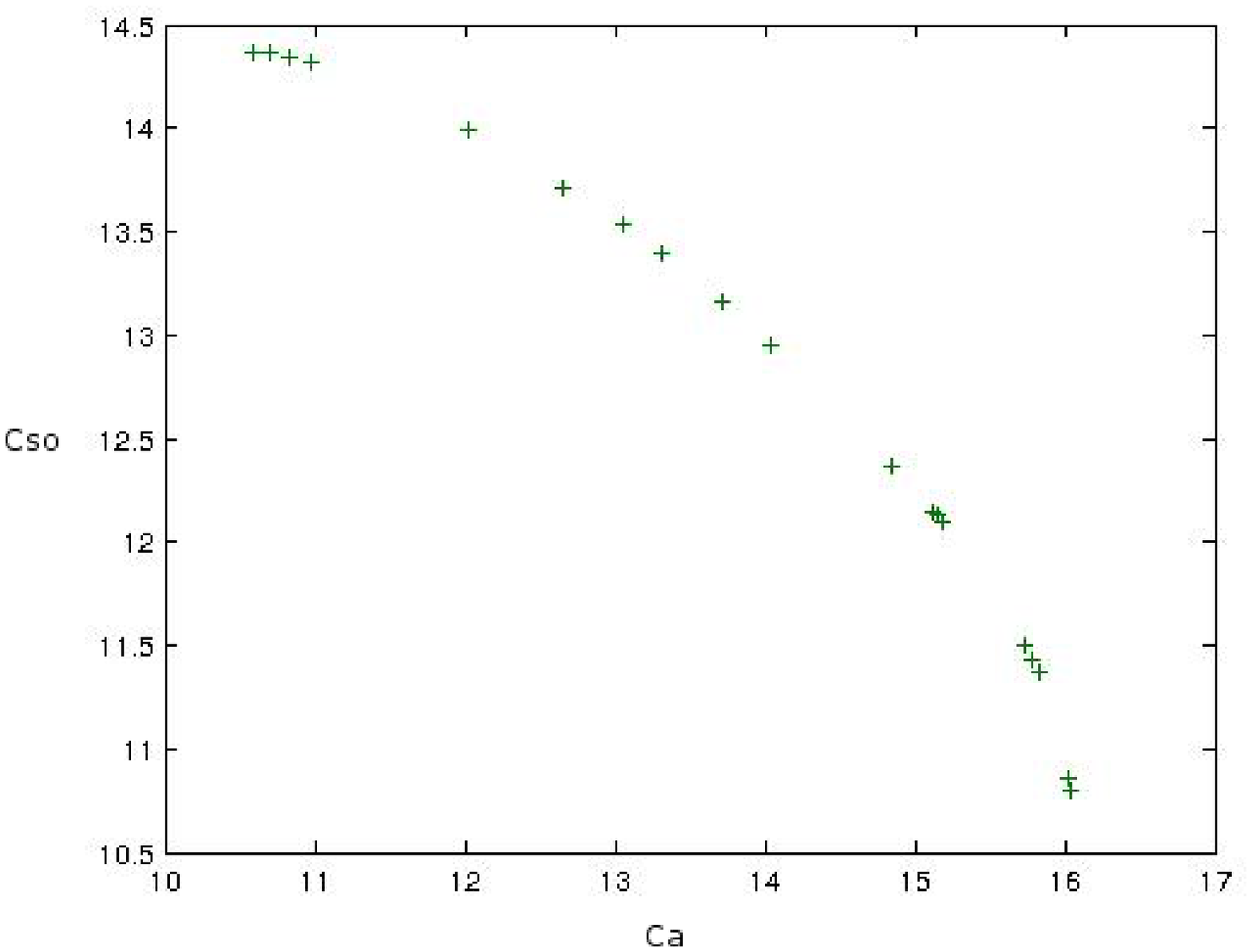}
	  \\
      (a) & (b) \\
    \end{tabular}
    \caption*{Fig 2. Weighted sum method (a) and reference point method (b) applied to the multi-species management problem using respectively 20 equally distributed pairs of weights and reference points. Ca is the sum over 20 years of the normalized density of abalone (divided by the maximal density). Cso is the sum over 20 years of the normalized number of sea otters (divided by the maximal number). \label{fig:A_SO}}
\end{figure}

Note that in both cases the computing time was very small and not reported here. This is because both cases were modelled by linear programs using only continuous variables, typically fast to solve~\cite[chap. 8]{bertsimas1997introduction}.



\subsection{Spatial allocation of resources}
\label{sub:resSAR}

We compared an explicit approach with a reference point approach.

An analysis of the 300 pairs of points corresponding to the explicit approach and the reference point approach revealed that the average minimum relative gap between an explicit point and its projection obtained using our reference point approach was $27.74\%$. This means that solution points generated by the reference point method were, on average, at least $27.74\%$ greater on every criteria than the points provided by the explicit method.

This result is not surprising, because the reference point generates only non-dominated points (see the guarantees of Section~\ref{sub:RP_method}), while the explicit method has a very low probability of generating a non-dominated point.

Table~\ref{tab:compare_explicit} 
illustrates the superiority of the reference point method compared to the explicit method. 
\begin{center}
\begin{table}[h!]
\begin{centering}
\begin{tabular}{ccccccc}
\hboldline
\textbf{Pair} & \textbf{Method} &	\textbf{Water travel time}	&	\textbf{Carbon} 	&	\textbf{Species 1} & 	\textbf{Species 2}  & 	\textbf{Species 3}
\tabularnewline
\hboldline
1	&	Explicit	&	1637	&	512	&	564 & 551 & 580	\\
1	&	\textbf{RP}	& \textbf{2719} & \textbf{847} &  \textbf{897} &  \textbf{884} &  \textbf{913}	\\
11	&	Explicit	&	1590 & 507 & 656 & 493 & 505 \\
11	&	\textbf{RP}	& \textbf{2611} & \textbf{842} &  \textbf{989} &  \textbf{830} &  \textbf{838}	\\
22	&	Explicit	&	1532 & 620 & 537 & 533 & 570 \\
22	&	\textbf{RP}	& \textbf{2557} & \textbf{944} &  \textbf{862} &  \textbf{861} &  \textbf{894}	\\
\hboldline
\end{tabular}
\par\end{centering}
\caption{Comparation between a sampling-based multi-objective explicit approach and the reference point method through a spatial resource allocation problem. Among all generated pairs of points, three randomly selected pairs are compared in the criteria space (pairs 1, 11 and 22). Units are not relevant in this table since the data was randomly generated. 
\label{tab:compare_explicit}
}
\end{table}
\par\end{center}

The total computation time for both methods was very low. For the reference point method, which is of course the slowest of the two methods, generating all the 300 points took only 84 seconds, i.e. 0.28 seconds per point on average. 

%

\section{Discussion}

Two main types of method for solving multi-objective problems exist in conservation: methods solving simplistic decision problems but using elaborate multi-objective decision-making processes, e.g. \cite{regan2009conservation} and Section \ref{sub:expl_app}, and optimization methods solving complex problems but using simplified and inaccurate decision-making process, e.g. \cite{chades2012setting} and Section \ref{sub:aggreg}. This paper considers a new approach for reconciling these two extreme types of approaches: the reference point method coupled with linear programming. The method can optimally solve multi-objective combinatorial problems while using an accurate interactive decision-making process. 


The theoretical features of the reference point method unlock a large range of important issues of multi-objective decision-making in conservation such as ethics, significance, transparency, convenience, interactivity and optimality (see Section~\ref{sec:method}). Additionally, the method avoids classic assumptions about the decision-maker's preferences. 
Results from the two conservation problems show that the method outperforms classic approaches by providing either better guidance for the decision-maker or better solutions on every criteria (Section~\ref{sec:results}).

The main caveat of the method is the need for an efficient linear programming formulation of the problem. 
The development of such formulations needs strong linear programming modelling techniques \cite{bertsimas1997introduction,williams2013model}.  However, in the particular case of multi-objective problems using a Markov decision process formalism, one can directly use our general formulation provided in Section~\ref{sub:methDPC}, Box 2.

As \cite{walters1978ecological,kareiva2014theory} and more recently \cite{beger2015integrating} emphasize, there is a real need to find good compromise solutions for multi-objective conservation and ecological problems in general. 
The approach could also be used to extend single-objective optimization techniques that tackle adaptive management problems~\cite{williams2011adaptive} and decision problems under partial observability ~\cite{chades2008stop}, where interactive methods seem particularly relevant. For adaptive management, only simple methods have been investigated to date, either based on the explicit approach methodology \cite[page 292]{walters1978ecological}, or based on the weighted sum method~\cite{probert2011managing,williams2011adaptive}. 
Recent approaches to find good compromise solutions between simplicity and optimality in conservation ~\cite{dujardin2015alpha} should also benefit from our approach.


With the increasing need to account for multiple objectives in conservation, the linear programming reference point approach should positively impact the way of solving multi-objective decision problems involving complex systems.

\pagebreak

\end{document}